\RequirePackage{fix-cm}
\documentclass[smallextended]{svjour3}       
\smartqed  
\usepackage{graphicx}
\usepackage{amssymb,amsfonts}
\usepackage{amscd, pdfsync, vpe}
\usepackage{amsmath}

\newtheorem{prop}[theorem]{Proposition}
\newtheorem{coro}[theorem]{Corollary}
\newtheorem{exa}[theorem]{Example}
\newtheorem{defi}[theorem]{Definition}
\newtheorem{rem}[theorem]{Remark}
\newcommand{\bte}{\begin{theorem}\quad  }
\newcommand{\ete}{\end{theorem} }
\newcommand{\bpr}{\begin{prop}\quad  }
\newcommand{\epr}{\end{prop} }
\newcommand{\ble}{\begin{lemma}\quad }
\newcommand{\ele}{\end{lemma}}
\newcommand{\bco}{\begin{coro}\quad }
\newcommand{\eco}{\end{coro} }
\newcommand{\bex}{\begin{exa}\quad \rm }
\newcommand{\eex}{\end{exa} }
\newcommand{\bde}{\begin{defi}\quad \rm }
\newcommand{\ede}{\end{defi} }
\newcommand{\brm}{\begin{rem} \quad \rm}
\newcommand{\erm}{\end{rem} }
\newcommand{\bdm}{\begin{displaymath} }
\newcommand{\edm}{\end{displaymath} }

\newcommand{\lb}{\label}
\newcommand{\lo}{\longrightarrow}

\begin{document}

\title{ Chain conditions on (Rees) congruences of $S$-acts
}


\author{Roghaieh Khosravi         \and
       Mohammad Roueentan 
}


\institute{R. Khosravi \at
              Department of Mathematics, Faculty of Sciences, Fasa University, Fasa, Iran. \\
              \email{khosravi@fasau.ac.ir}           
 \and
 M. Roueentan \at
 College of Engineering, Lamerd Higher Education Center,
Lamerd, Iran.
}
\date{Received: date / Accepted: date}

\maketitle

\begin{abstract}
In this paper, we introduce (Rees) artinian $S$-acts as those acts that satisfy the descending chain condition on their (Rees) congruences. Then, we show that artinian $S$-acts are those acts whose all factor acts are finitely (Rees) cogenerated. In addition, we continue the study of (Rees) noetherian $S$-acts, namely, those acts whose (Rees) congruences are finitely generated. Moreover, as a useful
tool for the investigation of chain conditions, we prove that the properties of being (Rees) noetherian and artinian are inherited in Rees short exact sequences. Next, we specifically consider the chain conditions on monoids. Finally, we prove that every right Rees artinian, commutative monoid is right
Rees noetherian.
\end{abstract}

\section{Introduction}

The properties of being noetherian and artinian are of fundamental importance in abstract algebra.  Noetherian and artinian rings and modules have been widely studied. See \cite{and,wis} for example. The study of right noetherian semigroups was initiated by Hotzel in \cite{hotz}. Right noetherian semigroups were further studied by Kozhukhov
in \cite{kozh}. According to \cite{coherency,gould}, a monoid $S$ is called weakly right noetherian
if every right ideal is finitely generated, and it is said to be right noetherian if every right congruence is finitely
generated. Later on, some connections between noetherian properties of a monoid and
coherency were investigated. Some fundamental properties of right
noetherian semigroups can be found in \cite{miller}.

According to \cite{davvaz}, an $S$-act is called noetherian if every ascending chain of its subacts is eventually stationary. Dually, it is called artinian if every descending chain of its subacts is eventually stationary. Moreover, the descending and ascending chain conditions on commutative monoids with zero elements as $S$-acts have been discussed. As defined in \cite{uniform}, an $S$-act is called noetherian (respectively, strongly noetherian) if it satisfies the ascending chain condition for subacts (respectively, congruences).

Quite similar to what we did in \cite{f.cog} for the introduction of finitely (Rees) cogeneraed $S$-acts, in this paper we
use the properties of being noetherian and Rees noetherian instead of being strongly noetherian and noetherian, respectively.
It is clear that an $S$-act $A_S$ is Rees noetherian if and only if each of its subacts
is finitely generated. The analogous equivalent condition in the artinian case is that each of its factor acts
is finitely cogenerated.
As the dual concepts of generating sets and finitely generated $S$-acts, the concepts of cogenerating classes of $S$-acts and finitely (Rees) cogenerated $S$-acts were investigated in \cite{f.cog}.

In this paper, we introduce and study artinian and Rees artinian $S$-acts as the dual notions of being noetherian and Rees noetherian, respectively.
In fact, the paper is an investigation of several basic
properties of (Rees) noetherian and artinian $S$-acts. In the reminder of this section, we provide the necessary background in the form of definitions and
statements of some theorems established in \cite{f.cog}, without proofs. In Section 2, we define an $S$- act to be artinian (noetherian) if it satisfies the descending (ascending) chain condition on its congruences. Moreover, we define a Rees artinian (noetherian) $S$-act as one that satisfies the descending (ascending) chain condition on its subacts. Then, we present some results on (Rees) noetherian and artinian $S$-acts.
 The main result studies
the behavior of these notions in Rees short exact sequences. Using these results, in Section 3, we investigate (Rees) noetherian and (Rees) artinian monoids. 
Throughout the paper, $S$ and $A_S$ are used to denote a monoid and a right $S$-act, respectively.
Recall that an equivalence relation $\rho$ on an $S$-act $A_S$ is said to be
a \textit{congruence} on $A_S$ if $a \rho a'$ implies $as~\rho ~a's$
for any $a, a'\in A_S$ and $s \in S$.
The set of all congruences on $A_S$ is denoted by $Con(A)$. Obviously, $ \nabla_A=A\times A, \Delta_A=\{(a,a)|~a\in A\} \in Con(A)$. It is known that the lattice $(Con(A_S),\subseteq,\vee,\cap)$ is a complete
lattice whose greatest and least elements are $\nabla_A$ and $\Delta_A$, respectively. The \textit{direct sum} of congruences $\sigma$ and $\delta$ is another congruence $\rho=\sigma\oplus\delta$, where $\rho=\sigma\vee\delta$ and $\sigma\cap\delta=\Delta_A$.  

An $S$-act is said to be \textit{simple} if it has no subacts other
than itself, and it is called \textit{$\theta$-simple} if it has no subacts other than
itself and the one-element subact $\Theta$. In contrast with module
theory, not every congruence on an $S$-act is associated with a subact, and
so it is normally the case that a simple (or $\theta$-simple) $S$-act has non-trivial congruences. For more information on $S$-acts, we refer the
reader to \cite{kilp2000}.


We recall from \cite{f.cog} some preliminaries related to finitely (Rees) cogenerated $S$-acts which will be needed in the characterization of (Rees) artinian $S$-acts.

 \bde (\cite{f.cog}) Let $A_S$ be an $S$-act. \begin{itemize}
  \item[\rm{(i)}] We say that $A_S$ is \textit{finitely cogenerated} if for every monomorphism $A \overset{f}{\longrightarrow} \prod_{i\in I} A_i$,
 $$ A \overset{f}{\longrightarrow} \prod_{i\in I} A_i \overset{\pi}{\longrightarrow}\prod_{j\in J} A_j$$ is also a monomorphism for some finite subset $J$ of $I$. Equivalently, for any family of congruences $\{\rho_i|~ i\in I\}$ on $A_S$, if $\cap_{i\in I}\rho_i=\Delta_A$, then $\cap_{j\in J}\rho_j=\Delta_A$ for some finite subset $J$ of $I$.
 \item[\rm{(ii)}] We call $A_S$ \textit{finitely Rees cogenerated} whenever for any family of Rees congruences $\{\rho_{B_i}|~ i\in I\}$ on $A_S$, if $\cap_{i\in I}\rho_{B_i}=\Delta_A$, then $\cap_{j\in J}\rho_{B_j}=\Delta_A$ for some finite subset $J$ of $I$.

\end{itemize}\ede

\bpr {\rm{(}}\cite{f.cog}{\rm{)}} \lb{le002} Let $A_S$ be an $S$-act, and $\theta\in Con(A)$. Then, $A/\theta$ is finitely (Rees) cogenerated if and only if for any family of (Rees) congruences $\{\rho_i|~ i\in I\}$ on $A_S$, if  $\cap_{i\in I}\rho_i=\theta$, then $\cap_{j\in J}\rho_j=\theta$ for some finite subset $J$ of $I$.
\epr




\bpr \lb{pr1.10} {\rm{(}}\cite{f.cog}{\rm{)}} Every finitely (Rees) cogenerated $S$-act is a finite coproduct of indecomposable $S$-acts.\epr

   \section{  					(Rees)	Artinian and  Noetherian $S$-acts}

 In this section, we provide the necessary notations
and establish several results on $S$-acts satisfying the ascending (or descending) chain condition on congruences and subacts.

 \bde  Let $A_S$ be an $S$-act. \begin{itemize}
        \item[\rm{(i)}] We call $A_S$ \textit{artinian {\rm{(}}noetherian{\rm{)}}} if $Con(A)$ satisfies the descending (ascending) chain condition.
       \item[\rm{(ii)}]  We say that $A_S$ is \textit{Rees artinian {\rm{(}}Rees noetherian{\rm{)}}} if it satisfies the descending (ascending) chain condition on its Rees congruences.
             \end{itemize} \ede
           The notions of \textit{right {\rm{(}}Rees{\rm{)}} artinian {\rm{(}}noetherian{\rm{)}} monoids} are applied for a monoid $S$ with these properties as a  right $S$-act.   It can be easily checked that $A_S$ is Rees artinian (Rees noetherian) if and only if it satisfies the descending (ascending) chain condition on its subacts.

          The following theorem  presents a characterization of noetherian $S$-acts. Its proof is similar to that of \cite[Proposition 2.1]{miller}.

\bte \lb{te4.01} For an $S$-act $A_S$, the following statements are equivalent.
 \begin{itemize}
 	\item[\rm{(i)}] The $S$-act $A_S$ is noetherian.
 	\item[\rm{(ii)}] Every congruence of $A_S$ is finitely generated.
 	 \item[\rm{(iii)}] Every non-empty subset of $Con(A)$ contains a maximal element.

 	\end{itemize} \ete
 Now, we characterize artinian $S$-acts.
 \bte \lb{te4.1} For an $S$-act $A_S$, the following statements are equivalent.
 \begin{itemize}
 	\item[\rm{(i)}] The $S$-act $A_S$ is artinian.
 	\item[\rm{(ii)}] Every factor act of $A_S$ is finitely cogenerated.
 	 \item[\rm{(iii)}] Every non-empty subset of $Con(A)$ contains a minimal element.
 	\end{itemize} \ete
 \begin{proof}
 (i) $\lo$ (iii) Let $\mathfrak{A}$ be a non-empty set of congruences of $A_S$, and
 suppose that $\mathfrak{A}$ does not have a minimal element. Then, for an element $\rho$ of $\mathfrak{A}$, $\rho$ is not minimal. Therefore, the
set $\{\theta\in \mathfrak{A}|~ \rho\supset \theta\}$ is not empty. Now, we can use the axiom of choice to find an infinite descending chain
$\rho\supset \theta_1\supset \theta_2\supset...$
 of congruences on $A_S$, which is a contradiction.

 	(iii) $\lo$ (ii). Using  Proposition \ref{le002}, let
 $\theta\in Con(A)$ and  $\{\rho_i|~ i\in I\}$ be a family of congruences on $A$ such that $\cap_{i\in I}\rho_i=\theta$. We show that $\cap_{j\in J}\rho_j=\theta$, for some finite subset $J$ of $I$. By our assumption, $$\{\cap_{j\in K}\rho_j| ~K\subseteq I\ \rm{is~ a~finite~set}\}$$
 	has a minimal element $\cap_{j\in J}\rho_j$, where $J$ is a finite subset of $I$. It can be easily checked that $\cap_{j\in J}\rho_j=\theta$.

 	(ii)$\lo$ (i). Suppose that $A_S$ has a descending chain
 $\theta_1\supset \theta_2\supset...$
 of congruences on $A_S$. Set $\theta=\cap_{i\in \mathbb{N}}\theta_i$. Then, since $A/\theta$ is finitely cogenerated, there exists a finite subset $\mathbb{F}$ of $\mathbb{N}$ with $\theta=\cap_{i\in \mathbb{F}}\theta_i$. Thus, for the maximum element $n$ of $ \mathbb{F}$, $ \theta_n= \theta_{n+j}$ for each $j\in  \mathbb{N}$.

 \end{proof}



 The proof of the following theorem is a special case of the proofs of the previous theorems, and this is why it is omitted.

 \bpr For an $S$-act $A_S$, the following statements are equivalent.
 \begin{itemize}
 	\item[\rm{(i)}] The $S$-act $A_S$ is Rees artinian {\rm{(}}noetherian{\rm{)}}.
 	\item[\rm{(ii)}] Every factor act {\rm{(}}subact{\rm{)}} of $A_S$ is finitely Rees cogenerated {\rm{(}}generated{\rm{)}}.
 	\item[\rm{(iii)}] Every non-empty set of subacts of $A_S$ contains a minimal {\rm{(}}maximal{\rm{)}} element.
 \end{itemize} \epr

 It is clear that every artinian (noetherian) $S$-act is Rees artinian (Rees noetherian). By \cite[Example 3.1]{gould}, the property of being Rees noetherian does not imply being noetherian. The following example shows that the property of being Rees artinian does not imply being artinian.
 \bex
  Let $S = (\mathbb{N}, \min)^\varepsilon = (\mathbb{N}, \min) \dot{\cup} \{\varepsilon\}$, where $\varepsilon$ denotes the externally adjoint identity. Then, $K_S = S\setminus \{\varepsilon\}$ is a right ideal of $S$. The subacts of $K_S$ are $1S\subseteq 2S\subseteq 3S\subseteq \ldots$. Hence, $K_S$ is Rees artinian. In \cite[Example 3.4]{f.cog}, it is shown that $K$ is not finitely cogenerated. Hence, $K$ is not artinian.

 \eex

   \ble \lb{le4.1}For a monoid $S$, the following statements are true. \begin{itemize}
   	\item[\rm{(i)}] Every subact of a {\rm{(}}Rees{\rm{)}} artinian  {\rm{(}}noetherian{\rm{)}} $S$-act is {\rm{(}}Rees{\rm{)}} artinian {\rm{(}}noetherian{\rm{)}}.
   	 \item[\rm{(ii)}] Every factor act of a {\rm{(}}Rees{\rm{)}} artinian  {\rm{(}}noetherian{\rm{)}}  $S$-act is {\rm{(}}Rees{\rm{)}} artinian {\rm{(}}noetherian{\rm{)}}.
   	
   	\end{itemize}\ele

   \begin{proof}  We prove the artinian case; one can use a similar method for the proof of the noetherian case.

   (i). Let $B$ be a subact of $A_S$, and assume that $A_S$ is artinian. Let $\rho_1\supseteq \rho_2\supseteq\ldots $ be a chain of congruences on $B$. Then, $\sigma_i=\rho_i\cup \Delta_A$ is a congruence on $A_S$, for each $i\in \mathbb{N}$. Since $A_S$ is artinian, $\sigma_n=\sigma_k$ for some $n\in \mathbb{N}$ and for each $k\geq n$. It can be easily checked that $\rho_n=\rho_k$, and $B$ is artinian.

   (ii). Let $B_S$ be a factor act of $A_S$, and assume that $A_S$ is artinian. By Theorem \ref{te4.1}, every factor act of $A_S$ is finitely cogenerated. So, every factor act of $B_S$ is finitely cogenerated, and $B$ is artinian.

   The statements for Rees artinian (noetherian) $S$-acts can be proved similarly.
\end{proof}
The following corollary is a straightforward result, and can be proved using (ii) of the previous lemma.
\bco  \lb{co6.1} Let $S$ be a monoid. Then, $S_S$ is {\rm{(}}Rees{\rm{)}} artinian {\rm{(}}noetherian{\rm{)}} if and only if every cyclic $S$-act is {\rm{(}}Rees{\rm{)}} artinian {\rm{(}}noetherian{\rm{)}}.\eco


 Let $f: A_S\lo B_S$ and $g: B_S\lo C_S$ be $S$-morphisms. Recall from \cite{Rees} that
  the sequence $ A \overset{f}{\lo} B \overset{g}{\lo} C$ is said to be a \textit{Rees short exact sequence} if $f$ is one-to-one, $g$ is onto, and $\ker{g}=\mathcal{K}_{\rm{Im}f}$, where $\mathcal{K}_{\rm{Im}f}=( f(A)\times f(A))\cup \Delta_B$. The following theorem discusses the behavior of the properties of being (Rees) artinian and noetherian for Rees short
exact sequences.

\bte \lb{te3.12} Let $  A\lo B\lo C$ be a Rees short exact sequence of $S$-acts. Then, $B_S$ is {\rm{(}}Rees{\rm{)}} artinian {\rm{(}}noetherian{\rm{)}} if and  only if both $A_S$ and $C_S$ are {\rm{(}}Rees{\rm{)}} artinian {\rm{(}}noetherian{\rm{)}}. \ete
  \begin{proof} The necessity is clear by Lemma \ref{le4.1}.

   \textit{Sufficiency.} Let  $\rho_1\supseteq \rho_2\supseteq\ldots$ be a descending chain of congruences on $B$. Then, $\sigma_i= \rho_i\cap(f(A)\times f(A))$ is a congruence on $f(A)$ for each $i\in \mathbb{N}$. Since $f(A)$ is artinian, $\sigma_n=\sigma_{n+1}=\cdots$ for some $n\in \mathbb{N}$. On the other hand, take $\varepsilon_i=\rho_i\vee \ker g$  and $\gamma_i=\{(g(a),g(b))|~(a,b)\in \varepsilon_i\} $ for each $i\in \mathbb{N}$. It is routine to check that $\gamma_i\in Con(C)$.  Moreover, $\gamma_1\supseteq \gamma_2\supseteq\ldots$. Since $C$ is artinian, $\gamma_m=\gamma_{m+1}=\cdots$. Let $k=\max\{m,n\}$. We show that $\rho_k=\rho_{k+1}=\cdots$. If $(b_1,b_2)\in \rho_k$, then $(b_1,b_2)\in \rho_k\vee \ker g=\varepsilon_k$, and so $(g(b_1),g(b_2))\in g(\varepsilon_k)=\gamma_k=\gamma_{k+1}$. There exists $(b'_1,b'_2)\in \varepsilon_{k
   +1}$ such that $(g(b_1),g(b_2))=(g(b'_1),g(b'_2))$. We obtain $b_1\ker g ~b'_1~ \varepsilon_{k   +1}~ b'_2 ~\ker g~ b_2$, which implies that $(b_1,b_2)\in \varepsilon_{k
   +1}$. Thus $\rho_k\subseteq \varepsilon_{k
   +1}$, showing that $\varepsilon_k=\varepsilon_{k+1}$. Moreover, since $\sigma_k=\sigma_{k+1}$, $\rho_k\cap (f(A)\times f(A))=\rho_{k+1}\cap (f(A)\times f(A))$. From $\ker g=( f(A)\times f(A))\cup \Delta_B$ we deduce that $\rho_k\cap \ker g=\rho_{k+1}\cap \ker g$. Now, we show that $\rho_k\subseteq \rho_{k+1}$. If $(b, b')\in  \rho_k$, then $\rho_k\subseteq \varepsilon_{k
   +1}$ implies that $(b, b')\in  \varepsilon_{k+1}$. Let $n$ be the least positive integer such that
   $$b\rho_{k+1}~d_1\ker g ~d_2~ \rho_{k +1}~ d_3 ~\ker g~ d_4~...~\rho_{k+1}~b,$$ where $d_1,\ldots,d_n\in B$. First, we claim that $n\leq 2$. Otherwise, if $n> 2$, then $(d_1,d_2),(d_3,d_4)\in \ker g$.  From $\ker g=( f(A)\times f(A))\cup \Delta_B$, four cases $d_1=d_2$ or $d_1, d_2\in f(A)$ and $d_3=d_4$ or $d_3, d_4\in f(A)$  are obtained, and in each case, one can reduce the length $n$, which contradicts the minimality of $n$. Thus, $n\leq 2$. If $n=1$, then it is clear that $(b, b')\in  \rho_{k   +1}$. If $n=2$, we obtain $b\rho_{k+1}~d_1\ker g ~d_2~ \rho_{k   +1}~b'$. Since $\rho_{k+1}\subseteq \rho_k$, $d_1\rho_{k}~b\rho_k ~b'~ \rho_{k }~d_2$, and $(d_1,d_2)\in \rho_k$. So, $(d_1,d_2)\in \rho_k\cap \ker g$. From $\rho_k\cap \ker g= \rho_{k+1}\cap \ker g$ we get $(d_1,d_2)\in \rho_{k+1}\cap \ker g$. Hence $b\rho_{k+1}~d_1\rho_{k+1} ~d_2~ \rho_{k   +1}~b'$, and so $(b,b')\in \rho_{k+1}$. Therefore, $\rho_k=\rho_{k+1}=\cdots$, and $B$ is artinian.

For the case of being Rees artinian, let  $\rho_{B_1}\supseteq \rho_{B_2}\supseteq \ldots$ be a descending chain of Rees congruences on $B$. Then, $\sigma_i= \rho_i\cap(f(A)\times f(A))=\rho_{B_i\cap f(A)}$ is a Rees congruence on $f(A)$ for each $i\in \mathbb{N}$.  On the other hand, if $\varepsilon_i=\rho_i\vee \ker g$, then $\gamma_i=\{(g(a),g(b))|~(a,b)\in \varepsilon_i\} =\rho_{g(B_i)}$ is a Rees congruence on $C$. Now, a reasoning similar to that of the previous proof can be used to show that $B$ is Rees artinian.

For the case of being (Rees) noetherian, it suffices to replace ``$\supseteq$" with ``$\subseteq$" in the previous proof.
    \end{proof}

    \ble \lb{le6.2} Let $A_S$ be an $S$-act, 
%
and $A_1\subseteq A_2\subseteq ... \subseteq A_n=A$. Then, $A$ is {\rm{(}}Rees{\rm{)}} artinian {\rm{(}}noetherian{\rm{)}} if and only if  $A_1$ and
the factor $S$-acts $A_{i+1}/A_i$  are  {\rm{(}}Rees{\rm{)}} artinian {\rm{(}}noetherian{\rm{)}} for all $1\leq i \leq n-1$.
\ele

\begin{proof}


Applying the Rees short exact sequence $$  A_i\lo A_{i+1}\lo  A_{i+1}/A_i$$ for each $0\leq i\leq n$, the result follows.\end{proof}

   We know from \cite{Rees} that for a monoid $S$ with zero, $  A_1\lo A_1\prod A_2\lo A_2$  and $  A_1\lo A_1\coprod A_2\lo A_2$ are Rees short exact sequences of $S$-acts. Note that in the category $S$-Act$_0$,  $A_1\coprod A_2=A_1\cup A_2$, where $ A_1\cap A_2= \Theta$, and each $A_i$ is called a $0$-direct summand of $A$.

   Now, using Theorem \ref{te3.12},
 we extend Lemma 2.9 of \cite{davvaz} from monoids to arbitrary $S$-acts in the following results.
 We discuss the way the properties of being (Rees) artinian and noetherian behave
on products and coproducts.

   \bpr \lb{pr5.1} For a monoid $S$, the following statements are true.\begin{itemize}
   	\item[\rm{(i)}] If  $ A=\prod _{i\in I}A_i$ {\rm{(}}$ A=\coprod _{i\in I}A_i${\rm{)}}  is    artinian {\rm{(}}noetherian{\rm{)}}, then $I$  is finite
and each $A_i$  is   artinian {\rm{(}}noetherian{\rm{)}}.

\item[\rm{(ii)}]   $ A_S$  is   {\rm{(}}Rees{\rm{)}}  artinian {\rm{(}}noetherian{\rm{)}} if and only if $A\coprod \Theta$  is  {\rm{(}}Rees{\rm{)}}  artinian {\rm{(}}noetherian{\rm{)}}.
\item[\rm{(iii)}] If  $ A_S$  is   {\rm{(}}Rees{\rm{)}}  artinian {\rm{(}}noetherian{\rm{)}}, then $\coprod _{i=1}^{i=n}A$  is  {\rm{(}}Rees{\rm{)}}  artinian {\rm{(}}noetherian{\rm{)}} for each $n\in \mathbb{N}$.

   In the following statements, assume that $S$ contains a zero.
     	\item[\rm{(iv)}] If $A_1,...,A_n$  are $S$-acts, then $ A=\prod _{i=1}^{i=n}A_i$  ($ A=\coprod _{i=1}^{i=n}A_i$ ) is (Rees) artinian (noetherian) if and only if each $A_i$, $1\leq i\leq n$, is (Rees) artinian (noetherian).
   	\item[\rm{(v)}] If a direct product $ A=\prod _{i\in I}A_i$ of $S$- acts is   Rees artinian {\rm{(}}noetherian{\rm{)}}, then $I$  is finite
and each $A_i$  is  Rees artinian {\rm{(}}noetherian{\rm{)}}.

\end{itemize} \epr
\begin{proof}
(i). The artinian case follows by the definition of being finitely cogenerated and Proposition \ref{pr1.10}. Let $ A=\prod _{i\in I}A_i$ be noetherian. If $A_i\neq \Theta$ for only finitely many $i\in I$. Suppose that there exists a finite subset $K$ of $I$ such that $A_i=\{\theta_i\}$ for each $i\in I\setminus K$. Since $\prod _{j\in \setminus K}A_j=\{(\theta_j)_{j\in I\setminus K}\}$ is a one-element $S$-act, we may assume without
loss of generality that $\Theta=A_{i_0}=\prod _{j\in I\setminus K}A_j$ for some $i_0\in I\setminus K$. Then $ A=(\prod _{i\in K'}A_j)$, where $K'=K\cup\{i_0\}$.  Therefore, $I$ may be considered as a finite set. Otherwise, assume to the contrary that
$J=\{i_n|~n\in \mathbb{N}\}$ is an infinite sequence of distinct elements of $I$ such that $A_{i_n}\neq \Theta$. Let $J_n=I\setminus\{i_k|~1\leq k\leq n\}$,  $B_n=\prod _{j\in J_n}A_j$, $\pi_n: A\rightarrow B_n$ be the natural epimorphism, and $\rho_n=\ker\pi_n$ for each $n\in \mathbb{N}$. It can be easily verified that $\rho_1\subseteq \rho_2\subseteq ...$. By our assumption, $\rho_m=\rho_{m+1}=\cdots$ for some $m\in \mathbb{N}$. Thus, it is routine to see that $A_{i_{m+1}}$ is the one-element $S$-act $\Theta$, a contradiction.

(ii) and (iii) are direct consequences of  Lemma \ref{le6.2}.

(iv).  Let  $ A=\prod _{i\in I}A_i$  be  Rees artinian. Assume to the contrary that $I$ is infinite. Suppose that
$J=\{i_n|~n\in \mathbb{N}\}$ is an infinite sequence of distinct elements of $I$. Take $ C_n=\prod _{i\in I}B_i$ such that $$ B_i= \left\{
\begin{array}{ll}A_i, &i\notin J~\rm{or} (i=i_k~\rm{for}~k\geq j)\\ \Theta, & (i=i_k~\rm{for}~k< j)\end{array} \right..$$
 Then, $A\supset C_1\supset C_2\supset...$  is a strictly descending chain of subacts of $A_S$,
a contradiction. For the case of being Rees noetherian, take $$ B_i= \left\{
\begin{array}{ll}\Theta, &i\notin J~\rm{or} (i=i_k~\rm{for}~k\geq j)\\ A_i, & (i=i_k~\rm{for}~k< j)\end{array} \right.,$$ and  the result follows using a proof similar to the previous one.

Part (v) is an immediate consequence of the previous theorem.
 \end{proof}

In the category $S$-Act$_0$, the coproducts of $\theta-$simple $S$-acts are called semi-simple. Moreover, one can show that over a monoid $S$ with zero, an $S$-act is semi-simple if and only if each of its subacts is a $0$-direct summand. Using \cite[Proposition 3.10]{f.cog} and Proposition \ref{pr5.1} we obtain the following result.
   \bco \lb{co2.9} Let $S$ be a monoid with zero. For a semi-simple $S$-act $A_S$, the following statements are equivalent.  \begin{itemize}
   	\item[\rm{(i)}] $A_S$ is Rees artinian. \item[\rm{(ii)}] $A_S$ is finitely Rees cogenerated.	\item[\rm{(iii)}] $A_S$ is finitely generated.
    \item[\rm{(iv)}] $A_S$ is Rees noetherian.  \end{itemize} \eco



Using \cite[Proposition 3.8]{f.cog} we obtain the following result.
\bpr \lb{pr3.50} Every right Rees artinian $S$-act has only finitely many maximal right subacts.\epr

By \cite[Lemma 11.6]{and}, for an
artinian module $M$ and a homomorphism $f : M \rightarrow M$, there exists a natural number
$n$ such that $\rm{Im} f^n + \ker f^n = M$. In the case of $S$-acts, we use the lattice of congruences to present the following analogous statement.
 \bpr \lb{pr5.20} Let $A_S$  be an $S$-act and $f \in
 \rm{End}(A_S)$. If $A_S$ is (Rees) artinian, then $\ker{f^n}\vee\mathcal{K}_{\rm{Im}f^n} = \nabla_A$ for some $n\in \mathbb{N}$. Moreover, $f$ is an
automorphism if and only if it is a monomorphism.\epr
\begin{proof}   It is obvious that
$$Imf\supseteq Imf^2 \supseteq \ldots .$$
Since $A_S$ is artinian, there exists
 $n \in \mathbb{N}$ such that $Imf^n =Imf^{2n}$.
Let $(a,a')\in \nabla_A$. Then $(f^n(a),f^n(a'))\in \mathcal{K}_{\rm{Im}f^{2n}}$, and so $f^n(a)=f^n(a')$ or $f^n(a)=f^{2n}(b)\neq f^n(a')=f^{2n}(b')\in \rm{Im}f^{2n}$.
If $f^n(a)=f^n(a')$, it is clear that $(a,a')\in \ker{f^n}$. Otherwise, $a~ \ker{f^n} ~ f^{n}(b) ~\mathcal{K}_{\rm{Im}f^n}~ f^{n}(b') ~\ker{f^n}~ a'$. Therefore, $(a,a')\in \ker{f^n}~\vee~\mathcal{K}_{\rm{Im}f^n}$, and the result follows.
To prove the second part note that if $f$ is a monomorphism, then $\ker f^n = \Delta_A$.
 So $\mathcal{K}_{\rm{Im}f^n}  = \nabla_A$, which implies that $f$ is an epimorphism.\end{proof}
As an immediate consequence of the second part of the previous proposition, we deduce that every Rees artinian $S$-act is cohopfian.

Now, using Proposition \ref{pr5.20}, we can state Fitting's Lemma for $S$-acts.
\bpr [Fitting's Lemma] Let $A_S$ be a noetherian and Rees artinian $S$-act. Then, for each $f \in
 \rm{End}(A_S)$ there exists $n\in \mathbb{N}$ such that $$\ker{f^n}\oplus\mathcal{K}_{\rm{Im}f^n} = \nabla_A.$$
 \epr
\begin{proof} Let $A_S$ be both noetherian and Rees artinian. Then, by Proposition \ref{pr5.20}, $\ker{f^m}\vee\mathcal{K}_{\rm{Im}f^m} = \nabla_A$ for some $m \in \mathbb{N}$ . Furthermore, by \cite[Proposition 3.6]{uniform},
$\ker{f^k}\cap\mathcal{K}_{\rm{Im}f^k} = \Delta_A$ for some $k \in \mathbb{N}$. Take $l=\max\{m,k\}$. Then, it is clear that $\ker{f^l}\oplus\mathcal{K}_{\rm{Im}f^l} = \nabla_A$. 

 \end{proof}

\section{					Right (Rees) Artinian and  Noetherian Monoids}
In this section, we focus on a monoid $S$ as a right $S$-act.
In \cite[Proposition 2.19]{miller}, right noetherian monoids are characterized
as monoids over which
 every cyclic $S$-act is finitely presented, or equivalently,
 every finitely generated $S$-act is finitely presented. By analogy with that result, we characterize right artinian monoids.
\bpr For a monoid $S$, the following statements are equivalent.\begin{itemize}
   	\item[\rm{(i)}]
 The monoid $S$ is right artinian.
	\item[\rm{(ii)}]  There exists a generator artinian $S$-act.
	\item[\rm{(iii)}]  Every finitely generated $S$-act is artinian.
	\item[\rm{(iv)}]  Every finitely generated $S$-act is finitely cogenerated.\end{itemize}\epr
\begin{proof}  The implications (i)$\Rightarrow$ (ii) and  (iii)$\Rightarrow$ (iv) are obviously true.

 (ii)$\Rightarrow$ (iii).  Suppose  that a generator $G$  is artinian, and that $B$ is finitely generated. Then there exists an epimorphism $f:\coprod_{i=1}^{i=n}G\rightarrow B$. Now, applying part (ii) of Proposition \ref{pr5.1} and Lemma \ref{le4.1}, we deduce that $B$ is artinian.

 (iv)$\Rightarrow$ (i).  By the assumption, every cyclic $S$-act is finitely cogenerated. So, by Corollary \ref{co6.1}, $S$ is artinian.\end{proof}

Lemma 2.17 of \cite{miller} shows that the property of being right noetherian for monoids is closed
under homomorphic images. The following result generalizes this fact to the behavior of chain conditions on $S$-acts under homomorphic images of monoids.

\bpr \lb{pr3.55} Let $f:S\lo T$ be a monoid homomorphism, and let $A$ be a $T$-act. If $A$ is a right {\rm{(}}Rees{\rm{)}} artinian {\rm{(}}noetherian{\rm{)}} $S$-act,
 then it is also a right {\rm{(}}Rees{\rm{)}} artinian {\rm{(}}noetherian{\rm{)}} $T$-act. If $f$ is an epimorphism, the converse is also true. \epr
 \begin{proof}  Let $f:S\lo T$ be a monoid homomorphism, and let $A$ be a right $T$-act. Clearly, $A$ is also a right $S$-act via $as=af(s)$ for any $a\in A,~s\in S$. Now, suppose that $A_S$ is a right artinian or noetherian $S$-act. If $\rho$ is a congruence on $A$ as $T$-act, it can be easily checked that it is a right congruence of $A$ as $S$-act. Now, the descending (or ascending) chain condition on the congruences of $A$ as $S$-act implies that it is also valid as $T$-act. Thus, $A$ is a right artinian (or noetherian) $T$-act.

 In the case of being Rees artinian (noetherian), we replace congruence with subact. The second part can be proved using a similar reasoning. \end{proof}

 By analogy with \cite[Lemma 2.17]{miller}, the following corollary is a direct consequence of the proposition above. It shows that
  all the artinian and noetherian properties are closed under quotients.

\bco If $\rho$ is a congruence on a right {\rm{(}}Rees{\rm{)}} artinian {\rm{(}}noetherian{\rm{)}} monoid $S$, then $S/\rho$ is a right {\rm{(}}Rees{\rm{)}} artinian {\rm{(}}noetherian{\rm{)}} monoid.\eco

It is well-known that every simple or $\theta$-simple $S$-act is both right Rees artinian and right Rees noetherian. The following result shows that over groups or 0-groups, the properties of being Rees artinian and Rees noetherian are essentially equivalent.
					
 \bpr \lb{pr5.44} Let $S$ be a group or 0-group. Then, an $S$-act $A_S$ is right Rees artinian if and only if it is right Rees noetherian.                                            \epr
 \begin{proof}   Suppose that $S$ is a group. By \cite[Proposition 1.5.34.]{kilp2000}, all right acts over $S$ are completely reducible. Moreover, by \cite[Corollary 3.11.]{f.cog},
every finitely Rees cogenerated, completely reducible
$S$-act is finitely generated. On the other hand, if  $A_S$  is  right Rees noetherian, again $A$ is finitely generated. So, we have the chain $A_1\subseteq A_2\subseteq ... \subseteq A_n=A$ such that  $A_1$ is simple and each
factor $S$-act $A_{i+1}/A_i$  is $\theta$-simple. Thus, the result follows by applying Lemma \ref{le6.2}.

Now, suppose that $S$ is a 0-group. Then, every $S$-act is semi-simple. The result follows by Corollary \ref{co2.9}.
 \end{proof}
In the following two results, we consider commutative monoids.
\bte Every right Rees artinian, commutative monoid is right Rees noetherian.\ete
\begin{proof}
Suppose that  $S_S$ is Rees artinian. Using \cite[Lemma 4.2]{hollow} and the fact that $S$ is commutative, we find that $S$ is a group or is local. If $S$ is a group, then it is clear that $S$ is right Rees noetherian. Otherwise, suppose that $M$ is the unique maximal right ideal of $S$.
Consider the chain
$S \supseteq M \supseteq M^2 \supseteq \ldots $.
Since $S$ is right Rees artinian, $M^n=M^{n+1}=\cdots$ for some $n\in \mathbb{N}$.
But, $S$ is commutative and Rees artinian, showing that it contains a proper minimal two-sided ideal, which is easily checked to be the unique minimum ideal denoted by $K(S)$. From the chain $S \supseteq M \supseteq M^2 \supseteq ...\supseteq M^n\supseteq K(S)$, by Lemma \ref{le6.2} it follows
that $B_i = M^{i-1}/M^i$ is a right Rees artinian $S$-act.
So, $B_i = M^{i-1}/M^i$ is a right Rees artinian $S/M^i$-act. Since
 $S/M$ is a factor of $S/M^i$, $B_i$ is a right Rees artinian $S/M$-act by Proposition \ref{pr3.55}. Thus, as $S/M$ is a 0-group, $B_i$ is a right Rees noetherian $S/M$-act by Proposition  \ref{pr5.44}, and it is a Rees noetherian $S$-act by Proposition \ref{pr3.55}. Therefore,
 Lemma \ref{le6.2} allows us to conclude that $A_S$ is right Rees noetherian.
\end{proof}

Recall that $S$ satisfies condition $A$ if every right $S$-act satisfies the ascending chain condition for cyclic subacts.
\bpr\lb{pr2.4}
 Every commutative, artinian semigroup $S$ satisfies condition $A$. 	

\epr

\begin{proof}
 Suppose that $S/\rho_1 \subseteq S/\rho_2 \subseteq ...$ is an ascending chain of cyclic acts, where for any $n\in \mathbb{N}$, $\rho_n $ is a congruence on $S$. Thanks to the commutativity of $S$, we observe that $ \rho_1 \supseteq  \rho_2 \supseteq  ...$, and the result follows by the assumption.
\end{proof}

Next, we describe Rees artinian monoids in a special case.
 \bpr Every right Rees artinian, left cancellative monoid is a group or a 0-group.\epr
\begin{proof} Let $a$ be a non-zero element of a right artinian, left cancellative monoid $S$. The chain of right ideals $aS \supseteq a^2S \supseteq a^2S \supseteq...$ stabilizes. So, we obtain $a^nS = a^{n+1}S$ for
some $n\in \mathbb{N}$ . If $a^n = a^{n+1}b$, then the fact that $S$ is left cancellative allows us to conclude that $ab = 1$.\end{proof}

To conclude this section, we consider monoid endomorphisms of artinian, projective $S$-acts.
\bpr\lb{pr2.3}
Suppose that $A_S$ is an artinian, projective right $S$-act. Then, $T$ satisfies the descending chain condition on principal right ideals, where $T=End(A_S)$.
\epr
\begin{proof}
		Suppose that ${f_1 T} \supseteq {f_2 T} \supseteq \ldots $ is a descending chain of principal right ideals of $T$, where ${f_n}\in T$ for every $n\in \mathbb{N}$. Then, for any $n\in \mathbb{N}$ there exists $g_n\in T$ such that $f_n g_n= f_{n+1}$. Hence, $f_n A \supseteq f_n g_n A \supseteq f_{n+1} A$ for every $n\in \mathbb{N}$. Therefore, $ {f_1 A} \supseteq {f_2 A} \supseteq\ldots $. Consequently, $ f_k A=f_{k+1} A $ for some $k\in \mathbb{N}$, by the assumption. We claim that $f_k T=f_{k+1} T $. Since $A_S$ is projective, there exists a homomorphism $h:A\longrightarrow A$ such that $f_k =f_{k+1}h$. Thus $f_k T=f_{k+1}T$, which completes the proof.
	\end{proof}

\def\bibname{REFERENCES}


\begin{thebibliography}{99}

\bibitem{and}  Anderson, F., and  Fuller, K.: Rings and Categories of Modules, Springer-Verlag,
New York, (1992)
\bibitem{coherency}  Dandan, Y.,  Gould, V.,  Hartmann, M.,  Rusuk, N.,  Zenab,  R.E.: Coherency and Construction for Monoids,  Quart. J. Math. 71, 1461–1488  (2020)
\bibitem{davvaz} Davvaz, B. ,  Nazemian, Z.: Chain conditions on commutative monoids, Semigroup Forum 100, 732–742  (2020)


  \bibitem{gould}   Gould, V.: Coherent monoids, J. Australian Math. Soc., 53, 166–182  (1992)

    \bibitem{hotz}  Hotzel, E.: On semigroups with maximal conditions. Semigroup Forum, 11, 337-362   (1975)

\bibitem{kilp2000}  Kilp, M.,   Knauer, U.,   Mikhalev, A.:  Monoids, Acts and Categories. W. de gruyter. Berlin, (2000)
  \bibitem{kozh}  Kozhukhov,  I.: On semigroups with minimal or maximal condition on left congruences. Semigroup Forum, 21, 337-350  (1980)


\bibitem{hollow} Khosravi, R.,   Roueentan, M.: Co-uniform and hollow S-acts over monoids, Commun. Korean Math. Soc., 2021
\bibitem{f.cog} Khosravi, R. , Liang, X.,   Roueentan, M.: Finitely cogenerated and
Cogenerating S-acts, https://arxiv.org/abs/2109.00344

\bibitem{miller} Miller, C., Rusuk, N.: Right Noetherian Semigroups, https://arxiv.org/abs/1811.08408v4

   \bibitem{wis} Wisbauer,  R.: Foundations of module and ring theory, Gordon and Breach Science Publishers, Reading (1991)




\bibitem{Rees} Chen, Y.,    Shum, K.P.:  Rees short exact sequences of $S$-systems,  Semigroup Forum, 65, 141-148  (2002)

    \bibitem{uniform} Roueentan,  M.,  Sedaghatjoo, M.: On uniform acts over semigroups,
 Semigroup Forum, 97, 229-243    (2018)


\end{thebibliography}
\end{document}